\documentclass[12pt]{article}
\usepackage{mathptmx}
\usepackage{graphicx}
\usepackage{amsmath}
\usepackage{amsfonts}
\newtheorem{thm}{Theorem}[section]

\newtheorem{rmk}{Remark}[section]

\newtheorem{Corollary}{Corollary}[section]

\newtheorem{tb}{Table}
\usepackage{epstopdf}

\title{ On Approximation by Kantorovich Exponential Sampling Operators}
\author{Shivam Bajpeyi  \thanks{Department of Mathematics, Visvesvaraya National Institute of Technology, Nagpur, Nagpur-440010, India. \newline E-mail: shivambajpai1010@gmail.com }
 \and
  A. Sathish Kumar \thanks{Department of Mathematics, Visvesvaraya National Institute of Technology, Nagpur, Nagpur-440010, India. \newline E-mail: mathsatish9@gmail.com}
  }
\date{}

\begin{document}
\maketitle
\bibliographystyle{plain}
\abstract{In this article, we extend our study of Kantorovich type exponential sampling operators introduced in \cite{own}. We derive the Voronovskaya type theorem and its quantitative estimates for these operators in terms of an appropriate K-functional. Further, we improve the order of approximation by using the convex type linear combinations of these operators. Subsequently, we prove the estimates concerning the order of convergence for these linear combinations. Finally, we give some examples of kernels along with the graphical representations.

\endabstract



\section{Introduction}
\label{intro}
The sampling theory has been finding its great applications mainly in approximation theory and signal processing. One of the most prominent result of the sampling theory is named after Wittaker-Kotelnikov-Shannon (WKS). The WKS sampling theorem \cite{wks} provides an exact reconstruction formula for the band-limited functions. Butzer \cite{butzer2} initiated the study of generalized sampling series which generalizes the WKS theorem for the functions which are not necessarily band-limited. These operators have great importance in developing the models for the reconstruction of signals. The theory has been developed by many authors in various aspects, see eg.,
\cite{butzer1,bardaro4,costa1,vinti1,bardarob,buta,butb,ries}.

In the past three decades, the exponential sampling has become an advantageous tool to deal with the problems arising in the wide areas of mathematics as well as physics, see eg. \cite{bertero,casasent,gori,ostrowsky}. Bartero, Pike \cite{bertero} and Gori \cite{gori} proposed the \textit{exponential sampling formula} to approximate the Mellin band-limited functions having exponentially spaced data. This formula is widely considered as the Mellin-version of the well known \textit{Shannon sampling theorem} \cite{butzer2}. But, the pioneering idea of mathematical study of exponential sampling formula is credited to Butzer and Jansche (see \cite{butzer5}) by using the theory of Mellin transform. The separate study of Mellin's theory has been initiated by Mamedov in \cite{mamedeo} and carried forward in \cite{butzer3,butzer5,butzer7}. Bardaro et al. generalized the above theory by replacing the $lin_{c}$ function in the exponential sampling formula by the more general kernel function in \cite{bardaro7}. It allows to approximate a continuous function by using its values at the points $ (e^{\frac{k}{w}}).$ But in practice, it is difficult to obtain the sample values at these nodes exactly. Thus, it is advantageous to replace the value $ f( e^{\frac{k}{w}})$ by the mean value of $ f(e^{x})$ in the interval $ \big[ \frac{k}{w}, \frac{k+1}{w} \big]$ for $ k \in \mathbb{Z},\  w>0.$ This led us to introduce and analyze the Kantorovich version of the generalized exponential sampling series in \cite{own}. The Kantorovich type generalization of operators is a significant subject in approximation theory as the Lebesgue integrable functions can be approximated by using Kantorovich type operators. In the last few decades, the Kantorovich modifications of several operators have been constructed and analyzed, eg.\cite{bardaro10,bardaro5,bardaro6,butzer2,vinti2,costa5,sir,gbs,SK,vijay2,PNA}. We also refer some of the recent developments related to the theory of exponential sampling, see \cite{bardaro7,bardaro1,bardaro11,comboexp,bardaro8,own}.

In view of Corollary 2 of \cite{bardaro7}, it is evident that the convergence of the operator $(S_{w}^{\phi}f)$ to $f$ is of order $\mathcal{O}(w^{-r})$ for $f \in C^{(r)}(\mathbb{R}^+)$ under the assumptions that the higher order moments of the kernel are null. The situation becomes comparatively downer for the Kantorovich version of the above operator which is introduced and analyzed in \cite{own}. In Theorem 2 of \cite{own}, we derived the following asymptotic formula under the assumption that the first order discrete moment vanishes. For $f \in C^{(2)}(\mathbb{R}^+)$ and $w>0,$ we have
$$ \lim_{w \rightarrow  \infty} w \big[ (I_{w}^{\chi}f)(x) - f(x) \big] = \frac{(\theta f)(x)}{2} \ , \hspace{0.3cm} x \in \mathbb{R}^+.$$
It is clear from the above formula that the convergence of the operator $(I_{w}^{\chi}f)$ to $f$ is of order $\mathcal{O}(w^{-1})$ for $f \in C^{(2)}(\mathbb{R}^+).$ In view of Theorem 3.1 of section 3, it is evident that one can not improve the order of convergence even if the existence of the higher derivatives for the function is assured and the higher order moments for the kernel are vanishing on $\mathbb{R}^+.$ Subsequently, it is hard to find the examples of the kernel function such that the higher order moments of the kernel are null. This motivates us to investigate further about the order of convergence of the family of operators $(I_{w}^{\chi}(f,.))_{w >0}.$
\par
The idea of considering the linear combination of the operators is mainly inspired from the pioneer study of Butzer \cite{lcbut} and by many authors in \cite{comb,lc1,lc2,comboexp}. In this paper, we implement this constructive approach to improve the order of convergence for the Kantorovich exponential sampling operators avoiding the constraint that the higher order moments for the kernels must vanish on $\mathbb{R}^+.$
\par
The paper is organized as follows. In section 3, we obtain the asymptotic formula and the quantitative estimates for the above family of operators in terms of Peetre's K-functional. The section 4 is devoted to analyse the approximation properties of the linear combinations of Kantorovich exponential sampling operators. We also prove the better order of convergence for these operators. Finally we have shown the approximation of different functions by the Kantorovich exponential sampling operators and its linear combination of operators. The error estimates are also provided.

\section{Preliminaries}
In what follows, we denote by $C^{(r)}(\mathbb{R}^+),\ r \in \mathbb{N}$ the space of all functions such that upto $r^{th}$ order derivatives are continuous and bounded on $\mathbb{R}^+$. Let $C(\mathbb{R}^+$) be the space of all continuous and bounded functions on $\mathbb{R}^+$. A function $f \in C(\mathbb{R}^+$) is called log-uniformly continuous on $\mathbb{R}^+$, if for any given $\epsilon > 0,$ there exists $\delta > 0$ such that $|f(u) -f(v)| < \epsilon$ whenever $| \log u - \log v| \leq \delta,$ for any $u, v \in \mathbb{R}^{+}.$ We denote the space of all   log-uniformly continuous and bounded functions defined on  $\mathbb{R}^{+}$ by $\mathcal{C}(\mathbb{R}^+).$ Similarly, $\mathcal{C}^{(r)}(\mathbb{R}^+)$ denotes the space of functions which are n-times Mellin continuously differentiable and $\theta^{r}f \in \mathcal{C}(\mathbb{R}^+).$ We consider $M(\mathbb{R}^{+})$ as the class of all Lebesgue measurable functions on $\mathbb{R}^+$ and $L^{\infty}(\mathbb{R}^{+}) $ as the space of all bounded functions on  $\mathbb{R}^{+}$ throughout this paper.

For $1 \leq p < +\infty$, let $ L^p(\mathbb{R}^+$) be the space of all the Lebesgue measurable and $p$-integrable functions defined on $\mathbb{R}^+$ equipped with the usual norm $\Vert f \Vert_p$. For $c \in \mathbb{R}$, we define the space
$$X_c = \{f : \mathbb{R}^+ \rightarrow \mathbb{C} : f(\cdot)(\cdot)^{c-1} \in  L^1(\mathbb{R}^+)\}$$ equipped with the norm
$$\Vert f \Vert_{X_c} = \Vert f(\cdot)(\cdot)^{c-1} \Vert_1 = \int_0^{+\infty} |f(u)|u^{c-1}du.$$
The Mellin transform of a function $f \in X_c$ is defined by
$$\hat{M}[f](s) := \int_0^{+\infty} u^{s-1}f(u)\ du \ , \,  \ (s = c + it, t \in \mathbb{R}).$$
A function $f \in X_{c} \cap C(\mathbb{R}^+), c \in  \mathbb{R}$ is called Mellin band-limited in the interval $[-\eta, \eta],$ if $\hat{M}[f](c+iw) = 0$ for all $|w| > \eta ,\ \eta \in \mathbb{R}^+.$
\par
Let $f : \mathbb{R}^+ \rightarrow \mathbb{C}$ and $c \in \mathbb{R}.$ Then, Mellin differential operator $\theta_c$ is defined by
$$\theta_cf(x) := xf'(x) + cf(x), \ \ \ \ \ x \in  \mathbb{R}^+ .$$
We consider $ \theta f(x) := \theta_{0} f(x) $ throughout this paper.
The Mellin differential operator of order $r \in \mathbb{N}$ is
defined by $\theta_c^r := \theta_c(\theta_c^{r-1}), \, \, \   \theta_c^1 := \theta_c,.$ The basic properties of the Mellin transform can be found in \cite{butzer3}.

\subsection{Kantorovich type Exponential Sampling Operators}
Let $ x \in \mathbb{R}^{+}$ and $ w >0.$ The Kantorovich version of the exponential sampling operators is defined by (\cite{own})
\begin{equation} \label{main}
(I_{w}^{\chi}f)(x)= \sum_{k= - \infty}^{+\infty} \chi(e^{-k} x^{w})\  w \int_{\frac{k}{w}}^{\frac{k+1}{w}} f(e^{u})\  du, \ \
\end{equation}
where $ f: \mathbb{R}^{+} \rightarrow \mathbb{R}$ is locally integrable such that the above series is convergent for every $ x \in \mathbb{R}^{+}$. It is clear that for $f \in L^{\infty}(\mathbb{R}^{+}),$ the above series is well defined for every $x \in \mathbb{R}^{+}.$ Let $ \chi :\mathbb{R}^{+} \rightarrow \mathbb{R}$ be the kernel function which is continuous on $\mathbb{R}^{+}$ such that it satisfies the following conditions: \\
(i) For every $ x \in \mathbb{R}^{+},$
$$  \sum_{k=- \infty}^{+\infty} \chi(e^{-k} x^{w}) =1.$$\\
(ii) For some $r >0,$ $M_{r}(\chi) <+ \infty$ and
$$ \lim_{\gamma \rightarrow + \infty} \sum_{|k-\log(u)|>  \gamma} |\chi(e^{-k} u)| \ |k- \log(u)|^{r}=0,$$ uniformly with respect to $ u \in \mathbb{R}^{+}.$
\begin{rmk}
The condition (ii) implies that there holds
$$\lim_{\eta \rightarrow + \infty} \sum_{|k-\log(u)|>  \eta} |\chi(e^{-k} u)| \ |k- \log(u)|^{j}=0 ,\ \ j=0,1, \cdots r-1.$$
\end{rmk}
We define the algebraic moments of order $\nu$ for the kernel function $ \chi$ as
$$ m_{\nu}(\chi,u):= \sum_{k= - \infty}^{+\infty}  \chi(e^{-k} u) (k- \log(u))^{\nu}, \hspace{0.5cm} \forall \ u \in \mathbb{R}^{+}.$$
Similarly, the absolute moment of order $\nu$ can be defined as
$$ M_{\nu}(\chi,u):= \sum_{k= - \infty}^{+\infty}  |\chi(e^{-k} u)| |k- \log(u)|^{\nu},  \hspace{0.5cm} \forall \ u \in \mathbb{R}^{+}.$$
We define $ \displaystyle M_{\nu}(\chi):= \sup_{u \in \mathbb{R}^{+}} M_{\nu}(\chi,u). $

\section{Approximation Results}

In this section, we establish some direct results e.g., Voronovskaya type asymptotic formula and its quantitative estimation for the Kantorovich exponential sampling operators (\ref{main}).

\begin{thm}\label{theorem1}
Let $\chi$ be the kernel function and $f \in C^{(r)}(\mathbb{R}^+).$ Then we have
$$ [(I_{w}^{\chi}f)(x) - f(x)]= \sum_{i=1}^{r} \frac{(\theta^{i}f)(x)}{(i+1)! \ w^{i}} \Bigg[ \sum_{j=1}^{i+1} \binom{i+1}{j} m_{i-j+1}(\chi,x) \Bigg]+ R_{r}(x),$$
where $R_{r}(x)$ represents an absolutely convergent series given by
$$ \sum_{k= - \infty}^{+\infty}  \chi(e^{-k} x^{w}) \  w \int_{\frac{k}{w}}^{\frac{k+1}{w}} h \left( \frac{e^u}{x} \right) (u- \log x)^r \ du.$$
Moreover, we have $R_{r}(x)= o(w^{-r})$ as $w \rightarrow \infty.$
\end{thm}

\noindent\bf{Proof.}\rm \
For $f \in C^{(r)}(\mathbb{R}^+),$ the Taylor's formula in terms of Mellin derivatives (\cite{bardaro7}) can be written as
$$ f(e^u) = f(x) + (\theta f)(x) (u- \log x) +
\frac{(\theta^{2}f)(x)}{2!} (u- \log x )^2 + \cdots
+\frac{(\theta^{r} f)(x)}{r!} (u- \log x)^r (x)+ h \left( \frac{e^u}{x} \right) (u- \log x)^r,$$ where $h : \mathbb{R}^+ \rightarrow \mathbb{R}$ is a bounded function such that $\displaystyle \lim_{t \rightarrow 1} h(t) = 0.$ In view of (\ref{main}) we obtain
\noindent $[(I_{w}^{\chi}f)(x) - f(x)]$
\begin{eqnarray*}
&=& \sum_{k= - \infty}^{+\infty} \chi(e^{-k} x^{w})\  w \int_{\frac{k}{w}}^{\frac{k+1}{w}} \left( \sum_{i=1}^{r}\frac{(\theta^{i} f)(x)}{i !} (u- \log x)^{i} + h \left( \frac{e^u}{x} \right) (u- \log x)^r \right)du \\
&=& \sum_{i=1}^{r}\frac{(\theta^{i} f)(x)}{i !} \sum_{k= - \infty}^{+\infty} \chi(e^{-k} x^{w})\  w \int_{\frac{k}{w}}^{\frac{k+1}{w}} (u- \log x)^{i} \ du + R_{r}(x),
\end{eqnarray*}
where $\displaystyle R_{r}(x):= \sum_{k= - \infty}^{+\infty} \chi(e^{-k} x^{w})\  w \int_{\frac{k}{w}}^{\frac{k+1}{w}} h \left( \frac{e^u}{x} \right) (u- \log x)^r du. $
For any fixed index $\eta,$ we obtain
$$ \sum_{k= - \infty}^{+\infty} \chi(e^{-k} x^{w})\  w \int_{\frac{k}{w}}^{\frac{k+1}{w}} \frac{(\theta^{\eta} f)(x)}{\eta !} (u- \log x)^{\eta} du= \frac{(\theta^{\eta} f)(x)}{w^{\eta}(\eta+1)!} \Bigg[ \sum_{j=1}^{\eta+1} \binom{\eta+1}{j} m_{\eta-j+1}(\chi)\Bigg].$$
Next we estimate the remainder term $R_{r}(x).$
\begin{eqnarray*}
|R_{r}(x)| & \leq & \sum_{k= - \infty}^{+\infty} | \chi(e^{-k} x^{w}) | \  \left| w \int_{\frac{k}{w}}^{\frac{k+1}{w}} h \left( \frac{e^u}{x} \right) (u- \log x)^r \ du \right|\\
& \leq & \left( \sum_{| \frac{k}{w} - \log x | < \gamma} + \sum_{| \frac{k}{w} - \log x | \geq \gamma} \right) | \chi(e^{-k} x^{w}) | \  \left| w \int_{\frac{k}{w}}^{\frac{k+1}{w}} h \left( \frac{e^u}{x} \right) (u- \log x)^r \ du \right|.
\end{eqnarray*}
Using condition (ii) and the fact that $h(x)$ is bounded  such that $\displaystyle \lim_{x \rightarrow 1} h(x) = 0,$ we see that
$$  w^{r} |R_{r}(x)| \leq \frac{\epsilon}{(r+1)} \left( \sum_{j=1}^{r+1} \binom{r+1}{j} M_{r-j+1}(\chi) + \| h \|_{\infty} \right).$$
This clearly shows that the series is absolutely convergent and $R_{r}(x)= o(w^{-r})$ as $w \rightarrow \infty$ for arbitrary $\epsilon >0.$
On combining all the estimates, we obtain the desired result.\\

As a consequence of the above theorem, we establish the following corollary.
\begin{Corollary}
Let $\chi$ be the kernel and $f \in C^{(1)}(\mathbb{R}^+).$ Then, we have the following asymptotic formula
$$\displaystyle \lim_{w \rightarrow \infty}\ w [(I_{w}^{\chi}f)(x) - f(x)]=  \frac{(\theta f)(x)}{2} \big( 1+ 2 m_{1}(\chi,x)\big) + o(w^{-1}).$$
Moreover, if $m_{j}(\chi,x)=0$ for $1 \leq j \leq r-1$ then for $f \in C^{(r)}(\mathbb{R}^+),$ we have
$$[(I_{w}^{\chi}f)(x) - f(x)]= \sum_{i=1}^{r} \frac{(\theta^{i}f)(x)}{(i+1)! \ w^{i}}+ \frac{(\theta^{r}f)(x)}{(r!) \ w^{r}} m_{r}(\chi,x) + o(w^{-r}).$$
\end{Corollary}

\subsection{Quantitative estimates}
Now we derive the quantitative estimates concerning the order of convergence for the family of operators $(I_{w}^{\chi}(f,.))_{w >0}$ by using the Peetre's- K functional (see \cite{bardaro9}). The Peetre's K-functional for $f \in \mathcal{C}(\mathbb{R}^+)$ is defined by
$$ \hat{K}(f,\epsilon, \mathcal{C}(\mathbb{R}^+),\mathcal{C}^{(1)}(\mathbb{R}^+)) := \inf \{ \|f-g\|_{\infty} + \epsilon \|\theta g \|_{\infty} : g \in \mathcal{C}^{(1)}(\mathbb{R}^+), \epsilon \geq 0 \}. $$

\begin{thm}\label{theorem2}
Let $\chi$ be the kernel function and $f \in \mathcal{C}^{(1)}(\mathbb{R}^+).$ Then the following estimate holds
$$ \left| [(I_{w}^{\chi}f)(x) - f(x)] - \frac{(\theta f)(x)}{2w} \big( 1+ 2 m_{1}(\chi,x)\big)  \right| \leq \frac{1+ 2 M_{1}(\chi)}{w} \hat{K} \left( f, \frac{1}{6w} \left(\frac{1+3 M_{1}(\chi)+3 M_{2}(\chi)}{1+2 M_{1}(\chi)} \right)  \right).$$
 \end{thm}

\noindent\bf{Proof.}\rm \
From the Mellin's Taylor formula for $r=1,$ we have
\begin{eqnarray*}
\left| [(I_{w}^{\chi}f)(x) - f(x)] - \frac{(\theta f)(x)}{2w} \big( 1+ 2 m_{1}(\chi,x)\big)  \right| &=& \left| \sum_{k= - \infty}^{+\infty} \chi(e^{-k} x^{w}) \   w \int_{\frac{k}{w}}^{\frac{k+1}{w}}  h \left(\frac{e^u}{x} \right)(u- \log x) du \right| \\
& \leq & \sum_{k= - \infty}^{+\infty} |\chi(e^{-k} x^{w})| \left| w \int_{\frac{k}{w}}^{\frac{k+1}{w}}  h \left(\frac{e^u}{x} \right)(u- \log x) du \right|
:=I.
\end{eqnarray*}
Now we substitute $\displaystyle R_{1}(f,x,u)= h \left(\frac{e^u}{x} \right)(u- \log x),\ u \in \left[\frac{k}{w}, \frac{k+1}{w}\right].$ Then, by using the estimate $ \displaystyle |R_{1}(f,x,u)| \leq 2 |u- \log x| \ \hat{K} \left( \theta f, \frac{|u-\log x|}{4} \right)$ (see \cite{bardaro9}), we can write
\begin{eqnarray*}
I& \leq & \sum_{k= - \infty}^{+\infty} | \chi(e^{-k} x^{w})| \  w \int_{\frac{k}{w}}^{\frac{k+1}{w}}  2 |u- \log x| \hat{K} \left( \theta f, \frac{|u-\log x|}{4} \right) du \\
& = & 2 \| \theta(f-g) \|_{\infty} \sum_{k= - \infty}^{+\infty} | \chi(e^{-k} x^{w})| w \int_{\frac{k}{w}}^{\frac{k+1}{w}} |u- \log x| \ du + \frac{\| \theta^{2}g \|_{\infty}}{2}  \sum_{k= - \infty}^{+\infty} | \chi(e^{-k} x^{w})| w \int_{\frac{k}{w}}^{\frac{k+1}{w}} |u- \log x|^2 du \\
& = & \| \theta(f-g) \|_{\infty} \left( \frac{1+2 M_{1}(\chi)}{w} \right) + \| \theta^{2}g \|_{\infty} \left( \frac{1+3 M_{1}(\chi)+3 M_{2}(\chi)}{6 w^2 } \right) \\
&=& \left( \frac{1+2 M_{1}(\chi)}{w} \right) \left( \| \theta(f-g) \|_{\infty} + \| \theta^{2}g \|_{\infty} \ \frac{1}{6w} \left(\frac{1+3 M_{1}(\chi)+3 M_{2}(\chi)}{1+2 M_{1}(\chi)} \right) \right).
\end{eqnarray*}
Taking infimum over $g \in \mathcal{C}^{2}(\mathbb{R}^+),$ we obtain $$ I\leq \frac{1+ 2 M_{1}(\chi)}{w} \hat{K} \left( f, \frac{1}{6w} \left(\frac{1+3 M_{1}(\chi)+3 M_{2}(\chi)}{1+2 M_{1}(\chi)} \right)  \right).$$

\begin{thm}\label{t3}
Let $f \in \mathcal{C}^{(r)}(\mathbb{R}^+)$ and  $\chi$ be the kernel such that $m_{j}(\chi,x)=0$ for $1 \leq j \leq r-1 .$ Then, we have
$$ \Bigg | [(I_{w}^{\chi}f)(x) - f(x)]- \sum_{i=1}^{r} \frac{(\theta^{i}f)(x)}{(i+1)! \ w^{i}} - \frac{(\theta^{r}f)(x)}{(r!) w^{r}} m_{r}(\chi,x) \Bigg | \leq
\frac{2 A}{w^{r}(r+1)!} \ \hat{K} \left( f, \frac{B}{2A (r+1)w} \right),$$
where $ A= (1+ (r+1)M_{r}(\chi))$ and $B= (1+ (r+2)M_{r+1}(\chi)).$
\end{thm}

\noindent\bf{Proof.}\rm \
From the Mellin's Taylor formula, we have

\noindent $ \displaystyle \Big | [(I_{w}^{\chi}f)(x) - f(x)]- \sum_{i=1}^{r} \frac{(\theta^{i}f)(x)}{(i+1)! \ w^{i}}- \frac{(\theta^{r}f)(x)}{(r!) w^{r}} m_{r}(\chi,x) \Big |$
\begin{eqnarray*}
& \leq & \sum_{k= - \infty}^{+\infty} |\chi(e^{-k} x^{w})| \Big| w \int_{\frac{k}{w}}^{\frac{k+1}{w}}  h \left(\frac{e^u}{x} \right)(u- \log x)^r du \big| :=J.
\end{eqnarray*}
Putting $\displaystyle R_{r}(f,x,u)= h \left(\frac{e^u}{x} \right)(u- \log x)^r,\ u \in \left[\frac{k}{w}, \frac{k+1}{w}\right]$ and using the estimate $\displaystyle |R_{r}(f,x,u)| \leq 2 \frac{| u- \log x|^r}{r!} \ \hat{K} \left( \theta^r f, \frac{|u-\log x|}{2(r+1)} \right),$ we obtain
\begin{eqnarray*}
J & \leq & \frac{2 \| \theta^{r}(f-g) \|_{\infty}}{r!} \sum_{k= - \infty}^{+\infty} |\chi(e^{-k} x^{w})|   w \int_{\frac{k}{w}}^{\frac{k+1}{w}} |u- \log x|^r du \\
&&+ \frac{\| \theta^{r+1}g \|_{\infty}}{(r+1)!} \sum_{k= - \infty}^{+\infty} |\chi(e^{-k} x^{w})| w \int_{\frac{k}{w}}^{\frac{k+1}{w}} |u- \log x|^{r+1} du \\
& \leq & \frac{2 \| \theta^{r}(f-g) \|_{\infty}}{w^{r}(r+1)!} (1+ (r+1)M_{r}(\chi)) + \frac{\| \theta^{r+1}g \|_{\infty}}{w^{r+1}(r+2)!} (1+(r+2)M_{r+1}(\chi))\\
& \leq & \frac{2 (1+ (r+1)M_{r}(\chi))}{w^{r}(r+1)!} \left( \| \theta^{r}(f-g) \|_{\infty} + \frac{1}{2(r+1)w} \frac{(1+(r+2)M_{r+1}(\chi))}{(1+ (r+1)M_{r}(\chi))}\| \theta^{r+1}g \|_{\infty}  \right).
\end{eqnarray*}
Now taking infimum over $g \in \mathcal{C}^{2}(\mathbb{R}^+),$ we get the desired estimate.

\section{Construction of linear combinations}
This section is devoted to the study of approximation properties of the linear combinations of the Kantorovich exponential sampling operators. Our central aim is to construct the appropriate linear combination of the operators $(I_{w}^{\chi}(f,.))_{w >0}$ to produce the better order of convergence in the asymptotic formula. Let $c_{i}, \ i=1,2,...,p$ be non-zero real numbers such that $\displaystyle\sum_{i=1}^{p} c_{i}=1.$ For $x \in \mathbb{R}^+$ and $w >0,$ we define the linear combination of the above operators as follows.
\begin{eqnarray} \label{combo}
(I_{w,p}^{\chi})(f,x) &=& \sum_{i=1}^{p} c_{i}  \sum_{k= - \infty}^{+\infty} \chi(e^{-k} x^{iw})\  (iw)  \int_{\frac{k}{i w}}^{\frac{k+1}{iw}} f(e^{u})\  du
=\sum_{i=1}^{p} c_{i} \ (I_{iw}^{\chi})(f,x).
\end{eqnarray}
Now we prove the asymptotic formula for the family of operators defined in (\ref{combo}).

\begin{thm}\label{t4}
Let $\chi$ be the kernel function and $f \in C^{(r)}(\mathbb{R}^+).$ Then, we have
$$ [(I_{w,p}^{\chi})(f,x)- f(x)]= \sum_{i=1}^{p} \sum_{j=1}^{r} \frac{(\theta^{j}f)(x)}{(j+1)! \ w^{j}} \bar{M}_{j}^{p}(\chi)+ o(w^{-r}),$$
where $ \displaystyle\bar{M}_{k}^{p}(\chi) := \sum_{i=1}^{p} \frac{c_{i}}{i^{k}} \left( \sum_{j=1}^{k+1} \binom{k+1}{j} m_{k-j+1}(\chi,x) \right).$
\end{thm}

\noindent\bf{Proof.}\rm \
From the condition $ \displaystyle\sum_{i=1}^{p} c_{i}=1,$ we can write
$$ [(I_{w,p}^{\chi})(f,x)- f(x)]= \sum_{i=1}^{p} c_{i} \sum_{k= - \infty}^{+\infty} \chi(e^{-k} x^{iw})\  (iw)  \int_{\frac{k}{i w}}^{\frac{k+1}{iw}} (f(e^{u}) - f(x))\  du. $$
Now using the $r^{th}$ order Mellin's Taylor formula, we obtain
\begin{eqnarray*}
[(I_{w,p}^{\chi})(f,x)- f(x)] &=& \sum_{i=1}^{p} c_{i} \sum_{k= - \infty}^{+\infty} \chi(e^{-k} x^{iw})\  (iw)  \int_{\frac{k}{i w}}^{\frac{k+1}{iw}} \left( \sum_{j=1}^{r}\frac{(\theta^{j} f)(x)}{j !} (u- \log x)^{j} + h \left( \frac{e^u}{x} \right) (u- \log x)^r \right)\ du \\
&=& \sum_{i=1}^{p} c_{i} \sum_{k= - \infty}^{+\infty} \chi(e^{-k} x^{iw})\  (iw)  \int_{\frac{k}{i w}}^{\frac{k+1}{iw}} \left( \sum_{j=1}^{r}\frac{(\theta^{j} f)(x)}{j !} (u- \log x)^{j} du \right)\\&& +
 \sum_{i=1}^{p} c_{i} \sum_{k= - \infty}^{+\infty} \chi(e^{-k} x^{iw})\  (iw)  \int_{\frac{k}{i w}}^{\frac{k+1}{iw}} h \left( \frac{e^u}{x} \right) (u- \log x)^r \ du \\
&=& \sum_{i=1}^{p} c_{i} \left( \sum_{j=1}^{r}\frac{(\theta^{j} f)(x)}{w^{j} (j+1)! \ i^{j}} \left( \sum_{\eta=1}^{j+1} \binom{j+1}{\eta} m_{j-\eta+1}(\chi,x)  \right) \right)+ I,
\end{eqnarray*}
where  $\displaystyle I= \sum_{i=1}^{p} c_{i} \sum_{k= - \infty}^{+\infty}\chi(e^{-k} x^{iw})\  (iw)  \int_{\frac{k}{i w}}^{\frac{k+1}{iw}} h \left( \frac{e^u}{x} \right) (u- \log x)^r \ du.$ Now proceeding along the lines of Theorem \ref{theorem1}, it follows that $I= o (w^{-r})$ as $w \rightarrow \infty.$ Hence, the proof is completed. \\

From the above theorem, we deduce the following Voronovskaja type asymptotic results.
\begin{Corollary}
For $f \in C^{(1)}(\mathbb{R}^+),$ we have
$$ [(I_{w,p}^{\chi})(f,x)- f(x)] = \sum_{i=1}^{p} \frac{c_{i}}{i} \left( \frac{(\theta f)(x)}{2w} (1+ 2 m_{1}(\chi,x)) \right)+ o (w^{-1}).$$ Furthermore, if $m_{j}(\chi,x)=0$ for $1 \leq j \leq r-1,$ we have
$$ [(I_{w,p}^{\chi})(f,x)- f(x)]= \sum_{i=1}^{p} c_{i} \left( \sum_{j=1}^{r} \frac{(\theta^{j}f)(x)}{(j+1)! w^{j}\  i^{j}} \right) + \sum_{i=1}^{p} \frac{c_{i}}{i^{r}} \left( \frac{(\theta^{r}f)(x)}{r! \ w^{r}} m_{r}(\chi,x) \right) + o (w^{-r}).$$
\end{Corollary}

\begin{Corollary}
Under the assumptions of Theorem \ref{t4} and if $\bar{M}_{k}^{p}(\chi) =0$ for $k=1,2, \cdots p-1,$ then for $f \in C^{(r)}(\mathbb{R}^+)$ with $r \geq p,$ we have
$$ \lim_{w \rightarrow \infty} w^{p}[(I_{w,p}^{\chi} f)(x) - f(x)] = \frac{(\theta^{p}f)(x)}{p+1 !} \bar{M}_{p}^{p}(\chi).$$
\end{Corollary}
It is important to remark that $\bar{M}_{k}^{p}(\chi)$ does not vanish in general, even if the higher order moments for the kernel are zero on $\mathbb{R}^+.$ So, in order to have $\bar{M}_{k}^{p}(\chi)=0$ for $k=0,1,2, \cdots p-1,$ we need to solve the following system $$ \sum_{i=1}^{p} c_{i}=1,\  \sum_{i=1}^{p} \frac{c_{i}}{i}=0,\  \sum_{i=1}^{p} \frac{c_{i}}{i^{2}}=0,\ \cdots \sum_{i=1}^{p} \frac{c_{i}}{i^{p-1}}=0. $$
The solution of the above system yields a linear combination which provides the convergence of order at least $p$ for functions $f \in C^{(p)}(\mathbb{R}^{+}).$ Let $f \in C^{(2)}(\mathbb{R}^{+})$ and $p=2,$ then we have
$$ c_{1} + c_{2}=1 \ \ \ \mbox{and} \ \ \  c_{1}+ \frac{c_{2}}{2}=0. $$ On solving, we obtain $c_{1}= -1$ and $c_{2}= 2.$ Then the linear combination is given by
\begin{equation}\label{combo1}
(I_{w,2}^{\chi})(f,x) = - (I_{w}^{\chi})(f,x)+ 2 (I_{2w}^{\chi})(f,x)\ , \ \ \ x \in \mathbb{R}^+.
\end{equation}
Now for $r \geq 2,$ the asymptotic formula for combination (\ref{combo1}) acquires the form
$$[(I_{w,2}^{\chi})(f,x) - f(x)]=  \sum_{j=2}^{r} \frac{(\theta^{j}f)(x)}{(j+1)! w^{j}} \bar{M}_{j}^{2}(\chi) + o(w^{-r}),$$
which provides the convergence of order at least 2 for $f \in C^{(2)}(\mathbb{R}^{+}).$ Similarly, if we take $p=3$ then for $f \in C^{(3)}(\mathbb{R}^{+}),$ we need to solve the following system:
\begin{eqnarray*}
c_{1}+c_{2}+c_{3}&=& 1 \\
c_{1}+ \frac{c_{2}}{2}+\frac{c_{3}}{3} &=& 0 \\
c_{1}+\frac{c_{2}}{4}+\frac{c_{3}}{9} &=& 0.
\end{eqnarray*}
The solution gives $c_{1}= \frac{1}{2},$ $c_{2}= -4$ and $c_{3}=\frac{9}{2}.$
Then the following linear combination ensures the order of convergence atleast 3 for $f \in C^{(3)}(\mathbb{R}^{+})$
\begin{equation} \label{combo2}
(I_{w,3}^{\chi})(f,x) = \frac{1}{2} (I_{w}^{\chi})(f,x)+ (-4) (I_{2w}^{\chi})(f,x) + \frac{9}{2} (I_{3w}^{\chi})(f,x), \ \ \ x \in \mathbb{R}^+.
\end{equation}
The corresponding asymptotic formula has the form
$$ [(I_{w,3}^{\chi})(f,x) - f(x)]=  \sum_{j=3}^{r} \frac{(\theta^{j}f)(x)}{(j+1)! w^{j}} \bar{M}_{j}^{3}(\chi) + o(w^{-r}).$$

\begin{thm} \label{t5}
Let $\chi$ be the kernel function and $f \in \mathcal{C}^{(1)}(\mathbb{R}^+).$ Then, we have
$$ \Bigg| I_{w,p}(\chi,x) - f(x) - \sum_{i=1}^{p} \frac{c_{i}}{i} \left( \frac{(\theta f)(x)}{2w} (1+ 2 m_{1}(\chi,x)) \right) \Bigg| \leq \frac{(1+2 M_{1}(\chi))}{w} \sum_{i=1}^{p} \frac{c_{i}}{i} \hat{K} \left( \theta f, \frac{A}{B} \frac{1}{6w} \right)\ ,$$ where $ \displaystyle A= \sum_{i=1}^{p} \frac{c_{i}}{i^2}(1+3M_{1}(\chi)+3M_{2}(\chi)) $ and $\displaystyle B= \sum_{i=1}^{p} \frac{c_{i}}{i} (1+2 M_{1}(\chi)).$
\end{thm}

\noindent\bf{Proof.}\rm \
We see that
\noindent $\displaystyle \Bigg| I_{w,p}(\chi,x) - f(x) - \sum_{i=1}^{p} \frac{c_{i}}{i} \left( \frac{(\theta f)(x)}{2w} (1+ 2 m_{1}(\chi,x)) \right)\Bigg|$
\begin{eqnarray*}
& \leq & \sum_{i=1}^{p} c_{i} \sum_{k= - \infty}^{+\infty} | \chi(e^{-k} x^{iw})|  \ \Big| (iw)  \int_{\frac{k}{i w}}^{\frac{k+1}{iw}}  h \left(\frac{e^u}{x} \right) (u- \log x) du \Big| :=J.
\end{eqnarray*}
In an analogous way to Theorem \ref{theorem2}, we put $\displaystyle R_{1}(f,x,u)= h \left(\frac{e^u}{x} \right)(u- \log x)^r,\ u \in \left[\frac{k}{iw}, \frac{k+1}{iw}\right]$ for $i \in \mathbb{N}$ and using the estimate $\displaystyle|R_{1}(f,x,u)| \leq 2 | u- \log x| \ \hat{K} \left( \theta^r f, \frac{|u-\log x|}{4} \right)$, we obtain
\begin{eqnarray*}
J& \leq & \sum_{i=1}^{p} c_{i} \sum_{k= - \infty}^{+\infty} | \chi(e^{-k} x^{iw})| \  (iw) \int_{\frac{k}{iw}}^{\frac{k+1}{iw}}  2 |(u- \log x)| \hat{K} \left( \theta f, \frac{|u-\log x|}{4} \right) du \\
&=& \sum_{i=1}^{p} c_{i} \left( 2 \| \theta(f-g) \|_{\infty} \sum_{k= - \infty}^{+\infty} | \chi(e^{-k} x^{iw})| (iw) \int_{\frac{k}{iw}}^{\frac{k+1}{iw}} |u- \log x| du \right)+ \\&&
 \sum_{i=1}^{p} c_{i} \left( \frac{\| \theta^{2}g \|_{\infty}}{2}  \sum_{k= - \infty}^{+\infty} | \chi(e^{-k} x^{iw})| (iw) \int_{\frac{k}{iw}}^{\frac{k+1}{iw}} |u- \log x|^2 du \right)\\
& = & \sum_{i=1}^{p} c_{i} \left( \| \theta(f-g) \|_{\infty} \left( \frac{1+2 M_{1}(\chi)}{iw} \right) + \| \theta^{2}g \|_{\infty} \left( \frac{1+3 M_{1}(\chi)+3 M_{2}(\chi)}{6 i^{2}w^2 } \right) \right).
\end{eqnarray*}
Now taking infimum over $g \in \mathcal{C}^{2}(\mathbb{R}^+),$ we get the desired result. \\

Now, we have the following corollary.
\begin{Corollary}
Let $\chi$ be the kernel function and $f \in \mathcal{C}^{(r)}(\mathbb{R}^+)$ with $m_{j}(\chi,x)=0$ for $1 \leq j \leq r-1.$ Then, we have the following estimate \\
\noindent  $  \Bigg| \displaystyle[I_{w,p}(\chi,x) - f(x)] - \sum_{i=1}^{p} c_{i} \left( \sum_{j=1}^{r} \frac{(\theta^{j}f)(x)}{(j+1)! w^{j}\  i^{j}} \right) - \sum_{i=1}^{p} \frac{c_{i}}{i^{r}} \left( \frac{(\theta^{r}f)(x)}{r! \ w^{r}} m_{r}(\chi,x) \right)\Bigg|$
\begin{eqnarray*}
\leq \frac{2D}{w^{r} (r+1)!}   \sum_{i=1}^{p} \frac{c_{i}}{i^{r}} \hat{K} \left( \theta^{n}f, \frac{1}{w(r+1)} \frac{E}{F} \right),
\end{eqnarray*}
where $\displaystyle D= ( 1+ (r+1) M_{r}(\chi)) ,$ \ $\displaystyle E= \sum_{i=1}^{p} \frac{c_{i}}{i^{r+1}} \left( 1+ (r+2) M_{r+1}(\chi)+ \frac{(r+1)(r+2)}{2} M_{r}(\chi) \right)$ and $\displaystyle F=\sum_{i=1}^{p} \frac{c_{i}}{i^{r}} \left( 1+ (r+1) M_{r}(\chi) \right).$
\end{Corollary}

\section{Examples and graphical representation}
We begin with the well-known Mellin's B-spline function(see \cite{bardaro7}). The Mellin's B-splines of order $n$ for $ x \in \mathbb{R}^{+}$ are defined as
$$\bar{B}_{n}(x):= \frac{1}{(n-1)!} \sum_{j=0}^{n} (-1)^{j} {n \choose j} \bigg( \frac{n}{2}+\log x-j \bigg)_{+}^{n+1},$$
where $(x)_{+} := \max \{x,0 \}.$ Since $ \bar{B}_{n}(x)$ is compactly supported and continuous on $\mathbb{R}^+,$ we have $ \bar{B}_{n}(x) \in X_{c}$ for any $c \in \mathbb{R}.$ The Mellin transformation of $\bar{B}_{n}$ is given by
\begin{eqnarray} \label{splinefourier}
\hat{M}[\bar{B}_{n}](c+it) = \displaystyle \bigg( \frac{\sin(\frac{t}{2})}{(\frac{t}{2})} \Bigg)^{n} \ \ , \hspace{0.5cm} t \neq 0.
\end{eqnarray}
To show that the above kernel satisfies the assumptions (i) and (ii), we use the following Mellin's Poisson summation formula
$$ (i)^{j} \sum_{k= - \infty}^{+\infty} \chi(e^{k} x) ( k-\log(u))^{j} = \sum_{k= - \infty}^{+\infty} \hat{M}^{(j)}[\chi](2k \pi i) \ x^{-2 k \pi i} \ \ \ \  ,\ \ \ \ \ \ \mbox{for} \ k \in \mathbb{Z}.$$
The Mellin B-spline satisfies all the assumptions of the presented theory (see \cite{bardaro7,comboexp,own}).
Consider the Mellin B-spline of order 2
\begin{equation*}
\bar{B}_{2}(x) =
     \begin{cases}
      {1 -  \log x ,} &\quad\text{} \ \ \ \ {1 < x < e}\\
 {1+  \log x,} &\quad\text{} \ \ { \text{$e^{-1} < x < 1$}}\\
 {0,} &\quad\text{} \ \ \ \ \ { \text{otherwise}.}\\
\end{cases}
\end{equation*}

Now we show the approximation of $ f(x)= 1-\cos(4 e^x),\  x \in [0.5, 1]$ by $(I_{iw}^{\chi}f)(x)$ and $(I_{w,3}^{\chi}f)(x)$ obtained in (\ref{combo2}) for $i=1,2,3$ and $w=15$ (Fig.1). It is evident from the graph that the combination (\ref{combo2}) provides the better rate of convergence.

\begin{figure}[h]
\centering
{\includegraphics[width=1.1 \textwidth]{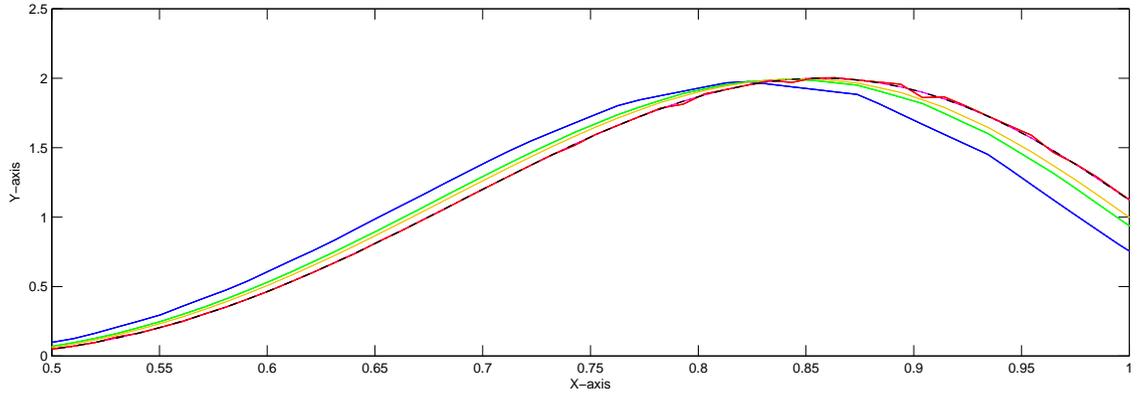}}
\caption{This figure shows the approximation of $ f(x)= (1-\cos(4e^x)), \ x \in [0.5, 1]$ (\textit{Black}) by $(I_{15}^{\bar{B}_{2}}f)(x)$ (Blue), $(I_{2\times15}^{\bar{B}_{2}}f)(x)$ (\textit{Green}), $(I_{3 \times 15}^{\bar{B}_{2}}f)(x)$ (\textit{Yellow}) and $(I_{15,3}^{\bar{B}_{2}}f)(x)$ (\textit{Red}) respectively.}
\end{figure}

\begin{tb}\label{table1}\centering
 {\it Error estimation (upto $4$ decimal points) in the approximation of $f(x)$ by $(I_{w,3}^{\chi}f)(x)$ and $(I_{iw}^{\chi}f)(x)$  for $i=1,2,3$ and $w=15.$}

$  $

\begin{tabular}{|l|l|l|l|l|}\hline
 $x$&$ |f(x) - I_{15}^{\bar{B}_{2}}f(x)|$&$|f(x)-I_{2\times15}^{\bar{B}_{2}}f(x)|$&$|f(x)-I_{3\times15}^{\bar{B}_{2}}f(x)|$&$|f(x)-I_{15,3}^{\bar{B}_{2}}f(x)|$\\
 \hline
 $0.60$&$0.1422$ & $0.0664$ & $0.0424$ & $0.0039$\\
  \hline
 $0.75$ & $0.1474$  &  $0.0807$ & $0.0561$ & $0.0033$\\
  \hline
 $0.80$ & $0.0613$  &  $0.0462$ & $0.0359$ & $0.0070$\\
 \hline
  $0.90$ & $0.2182$  &  $0.0800 $ & $0.0499$ & $0.0136$  \\  \hline
$0.95$ & $0.3230$  &  $0.1520 $ & $0.0963$ & $0.0129$  \\  \hline
             \end{tabular}
   \end{tb}

Now let us consider $4^{th}$ order Mellin B-spline function $\bar{B}_{4}(x).$ From the Mellin's Poisson summation formula, we obtain $m_{0}(\bar{B}_{4})=1,\  m_{1}(\bar{B}_{4})=0,\  m_{2}(\bar{B}_{4})= \frac{1}{3},\ m_{3}(\bar{B}_{4})=0.$  Now for $f \in C^{(2)}(\mathbb{R}^{+}),$ we have the asymptotic formula as
$$ \lim_{w \rightarrow + \infty} w[ (I_{w}^{\bar{B}_{4}}f)(x) - f(x)]= \frac{(\theta f)(x)}{2}.$$
But in view of (\ref{combo1}) and Corollary 2, we obtain the following asymptotic formula
$$\lim_{w \rightarrow + \infty} w^{2} [ I_{w,2}^{\bar{B}_{4}}f)(x) - f(x)]=  \frac{-(\theta^{2} f)(x)}{6}.$$
This clearly shows that the combination (\ref{combo1}) provides the order of convergence atleast 2 in asymptotic formula for $f \in C^{(2)}(\mathbb{R}^{+}).$
Subsequently, let $p=3$ and $f\in C^{(r)}(\mathbb{R}^+),\ r \geq 3.$ In view of Theorem 1, we have the following asymptotic formula
$$ \lim_{w \rightarrow + \infty} w[ (I_{w}^{\bar{B}_{4}}f)(x) - f(x)]= \frac{(\theta f)(x)}{2}.$$
But the combination (\ref{combo2}) provides the following asymptotic formula
$$ \lim_{w \rightarrow + \infty} w^{3} [ (I_{w}^{\bar{B}_{4}}f)(x) - f(x)]= \frac{(\theta^{3} f)(x)}{48} $$
which ensures the rate of convergence atleast 3 for $f \in C^{(3)}(\mathbb{R}^+).$
\begin{figure}[h]
\centering
{\includegraphics[width=1.1 \textwidth]{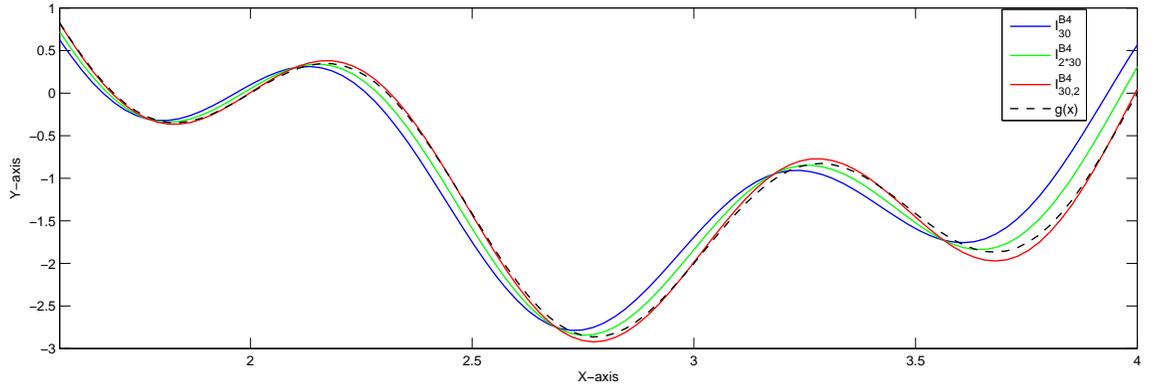}}
\caption{This figure shows the approximation of $ g(x)= \sin(2 \pi x)+2 \sin \left(\frac{\pi x}{2} \right),\  x \in \left[\frac{\pi}{2}, 4 \right]$ (\textit{Black}) by $(I_{30}^{\bar{B}_{4}}g)(x)$ (Blue), $(I_{2\times30}^{\bar{B}_{4}}g)(x)$ (\textit{Green}) and $(I_{30,2}^{\bar{B}_{4}}g)(x)$ (\textit{Red}) respectively.}
\end{figure}

\begin{tb}\label{table2}\centering
 {\it Error estimation (upto $4$ decimal points) in the approximation of $g(x)$ by $(I_{w,2}^{\bar{B}_{4}}f)(x)$ and $(I_{iw}^{\bar{B}_{4}}f)(x)$ for $i=1,2$ and $w=30.$}

$  $

\begin{tabular}{|l|l|l|l|}\hline
 $x$&$ |g(x) - I_{30}^{\bar{B}_{4}}g(x)|$&$|g(x)-I_{2\times30}^{\bar{B}_{4}}g(x)|$&$|g(x)-I_{30,2}^{\bar{B}_{4}}g(x)|$\\
 \hline
 $1.9$&$0.0880$ & $0.0385$ & $0.0110$\\
  \hline
 $2.6$ & $0.2217$  &  $0.1325$ & $0.0434$ \\
  \hline
 $3.1$ & $0.2037$  &  $0.1258$ & $0.0479$\\
 \hline
  $3.8$ & $0.4948$  &  $0.2071$ & $0.0806$  \\  \hline

             \end{tabular}
   \end{tb}

Now we consider the linear combination of the Mellin's B-spline functions of order $n$ as follows:
\begin{eqnarray*}
\chi(x) &:=& c_{1} [\tau_{\alpha} \bar{B}_{n}(x)] + c_{2}[ \tau_{\beta} \bar{B}_{n}(x)] \\
&=& c_{1}[\bar{B}_{n}(\alpha x)]+  c_{2}[\bar{B}_{n}(\beta x)] ,\hspace{0.3cm \forall x \in \mathbb{R}^+.}
\end{eqnarray*}
Here, $\tau_{h}$ represents h-translates of the function $\bar{B}_{n}(x)$ and is defined as $(\tau_{h}f)(x):= f(hx).$
The Mellin's transformation of $\chi(x)$ is given by
\begin{eqnarray} \nonumber
\hat{M}[\chi](w)&=& c_{1} \hat{M}[\bar{B}_{n}(\alpha x)]+ c_{2} \hat{M} [\bar{B}_{n}(\beta x)]\\
&=& c_{1} \ \alpha^{-w}  \hat{M} [\bar{B}_{n}](w)+ c_{2} \  \beta^{-w} \hat{M} [\bar{B}_{n}](w),
\end{eqnarray}
which gives
$$ \hat{M}[\chi](2k \pi i)= c_{1} \alpha^{- 2k \pi i}\hat{M}[\bar{B}](2k \pi i)+ c_{2} \ \beta^{-2k \pi i} \hat{M}[\bar{B}](2k \pi i).$$
Again from (5.6), we write
$$ \hat{M}^{'}[\chi](w)= c_{1} (\alpha^{-w} (\hat{M}^{'}[\bar{B}](w) - \alpha^{-w} \log \alpha) )+ c_{2} (\alpha^{-w} (\hat{M}^{'}[\bar{B}](w) - \beta^{-w} \log \beta) ).$$
In view of Lemma 1 in \cite{own}, we obtain
$$ c_{1} + c_{2}=1 \ , \ \ c_{1} \log \alpha + c_{2} \log \beta =0.$$
On solving for $c_{1} \ \ \mbox{and} \ \  c_{2},$ we get
$$ c_{1}= \frac{\log \beta}{(log \beta - \log \alpha)},\ \ c_{2}= \frac{- \log \alpha}{ (\log \beta - log \alpha)}.$$
In particular for $\alpha= e $ and $ \beta= e^{2},$ we obtain the following linear combination of the Mellin's B-spline of order 4
$$ \chi(x) = 2 \bar{B}_{4}(e x) - \bar{B}_{4}(e^{2} x)\ ,\ \ x \in \mathbb{R}^+.$$
Indeed, it satisfies all the assumptions of the kernel function. Again by using the Mellin's Poisson summation formula, we obtain
$$m_{0}(\chi)=1,\ m_{1}(\chi)=0,\ m_{2}(\chi)= \frac{-5}{3}.$$
Now for $f \in C^{(2)}(\mathbb{R}^+)$ and $p=2,$ we have the following asymptotic formulae for the operators $(I_{w}^{\chi}f) $ and $(I_{w,2}^{\chi}f),$ 
$$ \lim_{w \rightarrow + \infty} w[(I_{w}^{\chi}f)(x) - f(x)]= \frac{(\theta f)(x)}{2}$$ and
$$ \lim_{w \rightarrow + \infty} w^{2} [ (I_{w,2}^{\chi}f)(x) - f(x)]= \frac{1}{3}(\theta^{2} f)(x).$$
Evidently, the combination (\ref{combo1}) shows the convergence of order atleast 2 for $f \in C^{(2)}(\mathbb{R}^{+}).$


\end{document}